\documentclass[12pt,reqno]{amsart}
    \usepackage{amssymb,amsmath,amsthm,newlfont}
    \usepackage{psfrag}
    \usepackage[dvips,final]{graphicx}
\usepackage[T1]{fontenc}
\numberwithin{equation}{section}
\newcommand{\T}{\mathbb{T}}
\newcommand{\N}{\mathbb{N}}

\newcommand{\D}{\mathbb{D}}

\newcommand{\cJ}{{\mathcal{J}}}
\newcommand{\cI}{{\mathcal{I}}}

\newcommand{\cA}{{\mathcal{A}}}
\newcommand{\cH}{{\mathcal{H}}}
\newcommand{\cK}{{\mathcal{K}}}

\newcommand{\cL}{{\Lambda_\omega}}

\newtheorem{thm}{\bf Theorem}[section]
\newtheorem{coro}[thm]{\bf Corollary}

\newtheorem{prop}[thm]{\bf Proposition}

\newtheorem{rem}[thm]{\bf Remark}
\newtheorem{lem}[thm]{\bf Lemma}

\begin{document}

\title{{Closed ideals in analytic weighted Lipschitz algebras.} }

\date{\today}
\author[B. Bouya]{Brahim Bouya}
\address{Laboratoire Paul Painlev\'e, Universit\'e des Sciences et Technologies de Lille, B\^at. M2, 59655 Villeneuve d'Ascq Cedex, France.}
\email{bouya@math.univ-lille1.fr,\ Brahimbouya@gmail.com}

\thanks{This work was partially supported by the "Action Integr\'ee
Franco-Marocaine" No.\ MA/03/64.}

\begin{abstract}
{We obtain a complete description of closed ideals in weighted
Lipschitz algebras $\cL$ of analytic functions on the unit disk
satisfying the following condition
$$\frac{|f(z)-f(w)|}{\omega(|z-w|)}=o(1)\qquad(\mbox{as }|z-w|\longrightarrow0),$$
where $\omega$ is a modulus of continuity
satisfying some regularity conditions. In
particular the closed ideals of the algebras $\Lambda_{\chi_{_\alpha}},$
 where $\displaystyle \chi_\alpha(t):=\frac{1}{(|\log(t)|+1)^\alpha},$
$\alpha>0,$ are standard and this answers Shirokov's question \cite[page 587]{Shi1}.}
\end{abstract}
\maketitle

\section{\bf Introduction and statement of main result.}
Let $\D$ be the unit disk of the complex plane and $\T$ its
boundary. By $\cA(\D)$ we denote the usual disk algebra of all
analytic functions $f$ on $\D$ that are continuous on
$\overline{\D}.$ We define the weighted Lipschitz algebra
$\cL(\D)=\cL$ to be
$$\cL:=\Big\{f\in\cA(\D)\text{ :
}\sup_{z,w\in\D}\frac{|f(z)-f(w)|}{\omega(|z-w|)}=o(1)\quad(\mbox{as
}|z-w|\rightarrow0)\Big\}, $$ where $\omega(t)$ is a {\it modulus of
continuity}, i.e., a nondecreasing continuous real-valued function
on $[0,2]$ with $\omega(0)=0$ and  $\omega(t)/t$ is non
increasing function such that $\lim\limits_{t\rightarrow 0}\omega(t)/t=\infty.$ It is clear that $\cL$ is a commutative
Banach algebra when equipped with the norm
$$\|f\|_\omega:=\|f\|_\infty + \sup_{z,w\in\overline{\D}\atop z\neq w}\frac{|f(z)-f(w)|}{\omega(|z-w|)},$$
with $\|f\|_\infty:=\sup_{z\in\D}|f(z)|.$  Similarly the weighted
Lipschitz algebra $\Lambda_\omega(\T)$ is defined by
$$\Lambda_\omega(\T):=\Big\{f\in\cA(\D)\text{ :
}\sup_{
z,w\in\T}\frac{|f(z)-f(w)|}{\omega(|z-w|)}=o(1)\quad(\mbox{as
}|z-w|\rightarrow0)\Big\}.$$

Shirokov showed in \cite{Shi2} that $\cL$ possesses the so-called
F-property (Factorization property), i.e., for every given $f\in
\cL$ and inner function $U$ such that $f/U$ belongs to the algebra
$\mathcal{H}^{\infty}(\D)$ of bounded analytic functions, we have
$f/U\in \cL$ and $\|f/U\|_{\omega}\leq c \|f\|_{\omega},$
for an absolute constant $c$ \big(see
appendix B\big). Note that Tamrazov \cite{Tam}  proved that the
algebras $\cL$ and $\Lambda_\omega(\T)$ coincide for any arbitrary
modulus of continuity $\omega$ \big(see  appendix A\big).

The structure of closed ideals in the disk algebra is given independently
by Beurling and Rudin \cite{Hof}. They proved that if $\cI$ is a closed ideal of $\cA(\D),$ then there is an inner function
$U_{_\cI}$ (the greatest common divisor of the inner parts of the
non-zero functions in $\cI$) such that $\cI=\{f\in\cA(\D)\mbox{ :
} f_{|E_{_{\cI}}}\equiv0\mbox{ and
}f/U_{_\cI}\in\mathcal{H^{\infty}}(\D)\},$ where  $E_{_{\cI}}:=\{\xi\in\T\text{ :
}f(\xi)= 0 ,\ \forall f\in\cI\}.$

\noindent Using Beurling--Carleman--Domar resolvent's method
combined with the F-property,  we can reduce the problem of
characterization of closed ideals, in some algebras of analytic
functions, to a problem of approximation of outer functions (see for
example \cite{Bou} and references therein). Korenblum \cite{Kor} has
described the closed ideals of the algebra $H^2_1$ of analytic
functions $f$ such that $f'$ is in the Hardy space $H^2$. He proved
that these ideals are standard (in the sense of Beurling-Rudin
characterization of the closed ideals in the disk algebra). Later,
this result has been extended to some other Banach algebras of
analytic functions. In particular  by Matheson \cite{Mat} and
independently by Shamoyan \cite{Sha} for the algebras
$\Lambda_{\varphi_{_\alpha}},$ where
$\varphi_{_\alpha}(t):=t^\alpha,$ $0<\alpha<1.$\\
 The resolvent method is described
as follows: Define $d(\xi,E)$ to be the distance from $\xi\in\T$ to
the closed subset $E$ of $\T$ and let $\cI$ be a closed ideal of the
algebra $\Lambda_{\varphi_{_\alpha}}.$
\begin{itemize}
\item [1.] In the first step we
give an estimate to the norm of the resolvent $\|(\xi-\pi(z))^{-1}\|_{\Lambda_{\varphi_{_\alpha}}/\cI}$
in the quotient algebra $\Lambda_{\varphi_{_\alpha}}/\cI,$
where $\pi:\Lambda_{\varphi_{_\alpha}}\to \Lambda_{\varphi_{_\alpha}}/\cI$ is the canonical quotient map. We obtain
$\displaystyle \|(\xi-\pi(z))^{-1}\|_{\Lambda_{\varphi_{_\alpha}}/\cI}\leq \frac{c}{d^{4}(\xi,E_{_\cI})},$
 where $1\leq |\xi|\leq2$ and $c$ is an absolute constant.   So, from Cauchy formula on the quotient algebra $\Lambda_{\varphi_{_\alpha}}/\cI,$ we deduce that all functions in $\Lambda_{\varphi_{_\alpha}}$ such that $f/U_{_\cI}\in\cH^{\infty}(\D)$ and $\displaystyle |f(\xi)|\leq d^{4}(\xi,E_{_\cI}),$
$\xi\in\T,$ are in $\cI.$
\item [2.] The second step consists to prove that the space of all functions in
$\Lambda_{\varphi_{_\alpha}}$ such that
$f/U_{_\cI}\in\cH^{\infty}(\D)$ and $\displaystyle |f(\xi)|\leq
d^{4}(\xi,E_{_\cI}),$ $\xi\in\T,$ is dense in the standard ideal
$\{f\in\Lambda_{\varphi_{_\alpha}} \mbox{ : }
f_{|E_{_{\cI}}}\equiv0\mbox{ and
}f/U_{_\cI}\in\mathcal{H^{\infty}}(\D)\}.$
\end{itemize}
A closed subset $E\subset\T$ is called a Carleson set if the Carleson condition is satisfied, to wit
$$\frac{1}{2\pi}\int_0^{2\pi} \log\Big(\frac{1}{d(e^{it},E)}\Big)dt<+\infty.$$
 The zeros of each given function in any Banach algebras $H^2_1$, $\Lambda_{\varphi_{_\alpha}}$ and
 other ones in which the structure of closed ideals is also studied by using the
 resolvent method \cite{Bou, Shi1}, form a Carleson set. For the general case the resolvent method fails to apply,
 as example we can consider the algebras $\Lambda_{\chi_{_\alpha}},$  where
$\displaystyle \chi_\alpha(t):=\frac{1}{(|\log(t)|+1)^\alpha},$
$\alpha>0,$ \cite[page 587]{Shi1}. Indeed, let $\cI$ be a  closed
ideal of $\Lambda_{\chi_{_\alpha}}$ such that $E_{_\cI}$ is not a
Carleson set. We have $\displaystyle
\|(\xi-\pi(z))^{-1}\|_{\Lambda_{\chi_{_\alpha}}/\cI}\leq
\frac{c}{d^4(\xi,E_{_\cI})},$ where $1\leq |\xi|\leq2$ and
$\pi:\Lambda_{\chi_{_\alpha}}\to \Lambda_{\chi_{_\alpha}}/\cI$
is the canonical quotient map. It is clear that does not exist any
power $M$ such that $\displaystyle \int_\T
\frac{|f^M(e^{it})|}{d^4(e^{it},E_{_\cI})}dt<+\infty$ for all
functions in $\Lambda_{\chi_{_\alpha}}$ vanishing on  $E_{_\cI}$.
Therefore we can not conclude the first step of the resolvent
method as described above.

From now on,  $\omega$ will be a modulus of continuity such that for
every $1\leq\rho\leq2$ the following condition
\begin{equation}\label{cond}\omega(t^{\rho})\geq \eta_\rho\ \omega^{\rho}(t)\qquad(0\leq t\leq2),
\end{equation}
is satisfied, where  $\eta_\rho>0$ is a constant depending only on $\rho.$

In this work we prove that the closed ideals of the algebras $\cL$
are standard. For proving this we use only a special method of
approximating outer functions in $\cL$ together with the F-property.
More precisely, we obtain the following

\begin{thm}\label{prince} Let $\omega$ be a modulus of continuity satisfying \eqref{cond}.
If $\cI$ is closed ideal of $\cL$, then
$$\cI=\Big\{f \in\cL \text{ : } f_{|E_{_{\cI}}}\equiv0 \text{ and } f/U_{_{\cI}}\in\mathcal{H^{\infty}}(\D)\Big\},$$
where  $E_{_{\cI}}:=\{\xi\in\T\text{ : }f(\xi)= 0 ,\ \forall f\in\cI\}$
and $U_{_{\cI}}$ is the greatest common divisor of the inner parts
of the non-zero functions in $\cI.$
\end{thm}

Consequently, we obtain the structure of closed ideals of the
particular algebras $\Lambda_{\chi_{_\alpha}}.$

\section{\bf Other results and proof of Theorem \ref{prince}.}

We begin by recalling that every function $f$ in the disk algebra
has the canonical factorization $f=c_{_f}U_{_f} O_{_f},$  where
$c_{_f}$ is a constant of modulus 1, $U_{_f}$ is an inner function
(that is $|U_{_f}|=1$ a.e. on $\T$) and $O_{_f}$  the outer function
given by
$$O_{_f}(z)=\exp\Big\{\frac{1}{2\pi}\int_{0}^{2\pi}\frac{e^{i\theta}+z}{e^{i\theta}-z}\log|f(e^{i\theta})|d\theta\Big\}\qquad(z\in\D).$$
The closed ideal of all functions in $\cL$ vanishing on $E$ is
designed by $\cJ_{{E}}.$ In the section 3.2 we give the proof of the
following theorem

\begin{thm}\label{g2byf}
Let $\omega$ be a modulus of continuity satisfying the following condition
\begin{equation}\label{2condition}
\omega(t^{2})\geq \eta_{2}\ \omega^{2}(t)\qquad(0\leq t\leq2),
\end{equation}
where $\eta_{2}>0$ is a constant. Let $\cI$ be a closed ideal in $\cL$ such that $U_{_\cI}\equiv1$ and let $g\in\cJ_{E_{_\cI}}$ be an outer function.  Then $g^2$ belongs to $\cI.$
\end{thm}

\begin{rem}  In the same way, as in Theorem \ref{g2byf},  we can obtain that if moreover $\omega$  satisfy the following stronger condition
\begin{equation}\label{3condition}
\omega(t^{2})\geq \eta\ \omega(t)\qquad(0\leq t\leq2),
\end{equation}
then $g$ belongs to $\cI.$
\end{rem}

In the section 3.3 we prove the following theorem

\begin{thm}\label{gbyg2}Let $\omega$ be a modulus of continuity satisfying the condition \eqref{cond}.
Let $\cI$ be a closed ideal in $\cL$ and let $g\in\cL$ be a function such that $U_{_g}O^{2}_{_g}\in\cI.$ Then $g$ belongs to $\cI.$
\end{thm}

\textit{\bf Proof of Theorem \ref{prince}:} We have to prove that
every closed ideal of $\cL$ is standard. For this, let $\cI$ be a
closed ideal of the algebra $\cL.$ If $U_{_{\cI}}\equiv1$, let
 $g$  be a function in $\cJ_{E_{_\cI}}.$ Hence, making use of the F-property of $\cL$, it follows
$O_{_g}\in\cL,$ and therefore $O_{_g}\in\cJ_{E_{_\cI}}.$ Thus,
according to Theorem \ref{g2byf}, we deduce $O^2_{_g}\in\cI$ and
then $U_{_g}O^{4}_{_g}\in\cI.$ Next, by
applying Theorem \ref{gbyg2} two times, we conclude $g\in\cI.$\\
Now if  $U_{_{\cI}}\not\equiv 1,$ we choose $g\in \cJ_{E_{_\cI}}$
such that $g/U_{_\cI}\in\cH^{\infty}(\D).$ Thus, the associated
ideal
$$\cK_{_g}:=\{f\in\cL\ :\  fg\in\cI\}$$ is closed and by the
F-property of $\cL$ we  see easily that $U_{_{\cK_{_g}}}\equiv 1$
and then $\cK_{_g}=\cJ_{E_{_{\cK_{_g}}}}.$ Now, since
$E_{_{\cK_{_g}}}\subseteq E_{_\cI},$ then $O_{_g}\in\cK_{_g}.$
It follows that $U_{_g}O^2_{_g}\in\cI.$ So, by Theorem \ref{gbyg2},
$g\in\cI.$  This completes the proof of the theorem.

\section{\bf Approximation of functions in $\cL.$}

In this section we give the proof of Theorem \ref{g2byf} and  Theorem \ref{gbyg2}.
For simplicity we use the following Tamrazov's Theorem (see Appendix A): If $f$ is a function in the disk algebra such that $f\in\Lambda_{\omega}(\T),$ then $f\in\cL.$ We need also the following simple Lemma

\begin{lem}\label{shilbert}
Let  $f_n\in\cL$ be a sequence of functions converging uniformly on
the closed unit disk to $f\in\cL.$ If
$$\frac{|f_n(z)-f_n(w)|}{\omega(|z-w|)}=o(1)\qquad(as\ |z-w|\longrightarrow 0),$$
uniformly with respect to $n,$ then $\lim\limits_{n\longrightarrow
+\infty}\|f_n-f\|_{\omega}=0.$
\end{lem}

\subsection{Several technical lemmas}

For $f\in\cL,$ the inner function $U_{_f}$ is uniquely factored in the form $U_{_f}=B_{_f}S_{_f},$
 where $B_{_f}$ is the usual Blashke product associated to $Z_{_f}\cap\D,$  $Z_{_f}:=\{z\in\overline{\D}\text{ :
}f(z)=0\}$ and the function
$$S_{_f}(z):=\exp\Big\{-\frac{1}{2\pi}\int_{0}^{2\pi}\frac{e^{i\theta}+z}{e^{i\theta}-z}
d\mu_{_f}(\theta)\Big\}\qquad(z\in\D),$$
is the singular inner function associated to the singular positive
measure $\mu_{_f}.$ Note that the support $\mathrm{supp}(\mu_{_f})$
of the singular measure $\mu_{_f}$  is a closed subset of $E_{_f}:=\{\xi\in\T\ :\ f(\xi)=0\}.$ For $a,b\in\T,$ we design by $(a,b)$ (resp. $[a,b]$) an open arc
(resp. closed arc) of $\T$ connecting the points $a$ and $b.$

\begin{lem} \label{lem1}
Let $\omega$ be a modulus of continuity satisfying the condition \eqref{2condition}. Let $g$ be a function in $\cL$ and let $U=B_{_U}S_{_U}\in\cH^{\infty}(\D)$ be an inner function such that $B_{_g}/B_{_U}\in\cH^{\infty}(\D),$ $\mathrm{supp}(\mu_{_U})\subseteq E_{_g}$ and $(1/2\pi)\int_{0}^{2\pi}d\mu_{_U}(\theta)\leq M,$ where $M$ is a constant. Then $UO_g^2$ belongs to $\cL$ and we have
\begin{equation}\label{bb1}
\frac{|U(\xi)O_g^2(\xi)-U(\zeta)O_g^2(\zeta)|}{\omega(|\xi-\zeta|)}=o(1)\qquad (\mbox{as } |\xi-\zeta|\longrightarrow 0),
\end{equation}
uniformly with respect to $U.$
If moreover $\omega$ satisfy the condition \eqref{3condition}, then $UO_g$ belongs to $\cL$ and we have
\begin{equation}\label{bb2}
\frac{|U(\xi)O_g(\xi)-U(\zeta)O_g(\zeta)|}{\omega(|\xi-\zeta|)}=o(1)\qquad (\mbox{as } |\xi-\zeta|\longrightarrow 0),
\end{equation}
uniformly with respect to $U.$
\end{lem}

\begin{proof}
Let $\xi,\zeta\in\T$ be two distinct points such that $d(\xi,E_{_g})\geq d(\zeta,E_{_g}).$  It is clear that
\begin{eqnarray*}
&&\frac{|U(\xi)O_g^2(\xi)-U(\zeta)O_g^2(\zeta)|}{\omega(|\xi-\zeta|)}\\&\leq&
\frac{|B_{_U}(\xi)O_g^2(\xi)-B_{_U}(\zeta)O_g^2(\zeta)|}{\omega(|\xi-\zeta|)}+ |O_g^2(\zeta)|\frac{|S_{_U}(\xi)-S_{_U}(\zeta)|}{\omega(|\xi-\zeta|)}.
\end{eqnarray*}
By the F-property of $\cL$ we have $O_g\in\cL$ and  $B_{_U}O_g^2\in\cL.$ Then
to prove \eqref{bb1}, it is sufficient to prove that
\begin{equation}\label{singp}
|O_g^2(\zeta)|\frac{|S_{_U}(\xi)-S_{_U}(\zeta)|}{\omega(|\xi-\zeta|)}=o(1)\qquad(\mbox{as } |\xi-\zeta|\longrightarrow0).
\end{equation}
First we suppose that $\displaystyle|\xi-\zeta|\geq \big(\frac{d(\zeta,E_{_g})}{2}\big)^2.$ Then
\begin{eqnarray}\nonumber
|O_g^2(\zeta)|\frac{|S_{_U}(\xi)-S_{_U}(\zeta)|}{\omega(|\xi-\zeta|)}&\leq&
8\eta_{2}^{-1}\Big(\frac{|O_g(\zeta)|}{\omega(d(\zeta,E_{_g}))}\Big)^2 \\\label{sing1}
&=&o(1)\qquad (\mbox{as } |\xi-\zeta|\longrightarrow 0).
\end{eqnarray}
Now, let suppose that $\displaystyle|\xi-\zeta|\leq \big(\frac{d(\zeta,E_{_g})}{2}\big)^2.$ Then $[\xi,\zeta]\subset\T\setminus E_{_g}$ and therefore $d(z,E_{_g})\geq d(\zeta,E_{_g}),$ for every
$z\in[\xi,\zeta].$ There is $z\in[\xi,\zeta]$ such that
$$ \frac{|S_{_U}(\xi)-S_{_U}(\zeta)|}{|\xi-\zeta|}=|S_{_U}^{'}(z)|\leq \frac{1}{\pi}\int_{0}^{2\pi}\frac{1}{|e^{i\theta}-z|^2}d\mu_{_U}(\theta)\leq \frac{2M}{d^2(z,E_{_g})}.$$
It follows  $$\frac{|S_{_U}(\xi)-S_{_U}(\zeta)|}{|\xi-\zeta|}\leq \frac{2M}{d^2(\zeta,E_{_g})}.$$
We obtain
\begin{eqnarray}\nonumber
&& |O_g^2(\zeta)|\frac{|S_{_U}(\xi)-S_{_U}(\zeta)|}{\omega(|\xi-\zeta|)}\\\nonumber&=&
|O_g^2(\zeta)|\frac{|S_{_U}(\xi)-S_{_U}(\zeta)|}{|\xi-\zeta|}\frac{|\xi-\zeta|}{\omega(|\xi-\zeta|)} \\\nonumber&\leq&
2M\eta_{2}^{-1} \Big(\frac{|O_g(\zeta)|}{\omega(d(\zeta,E_{_g}))}\Big)^2 \frac{\omega(d^2(\zeta,E_{_g}))}{d^2(\zeta,E_{_g})}\frac{|\xi-\zeta|}{\omega(|\xi-\zeta|)}\\\label{sing2}&=&
o(1)\qquad(\mbox{as } |\xi-\zeta|\longrightarrow 0).
\end{eqnarray}
So \eqref{singp} follows from \eqref{sing1} and \eqref{sing2}. Consequently  $UO_g^2$ belongs to $\cL.$ If moreover $\omega$ satisfy the condition \eqref{3condition}, we can argue similarly to prove \eqref{bb2}. This finishes the proof of the lemma.
\end{proof}

\begin{lem}\label{lem3}
Let $f$ be a function in $\cL.$  Let  $\delta>0,$ $N\in\N$ and
$\{a_n\ :\ 0\leq n\leq N\}$ be a finite number of points in $E_{_f}.$ Then $$\lim\limits_{\delta\longrightarrow 0}\|\psi_{\delta,N}f-f\|_{\omega}=0,$$ where
$\displaystyle \psi_{\delta,N}(z):=\prod_{n=0}^{n=N}\frac{z\overline{a_n}-1}{z\overline{a_n}-1-\delta},\ z\in\D.$
\end{lem}

\begin{proof}
Without loss of the generality we can suppose that $N=0$ and
$a_0=1.$ Set $\displaystyle \psi_{\delta}(z):=\frac{z-1}{z-1-\delta},\ z\in\D,$
and suppose that $\{1\}\in E_{_f}.$ We have to show that
$\lim\limits_{\delta\longrightarrow
0}\|\psi_{\delta}f-f\|_{\omega}=0.$ Let $\xi,\zeta\in\T$ be two distinct points such that  $|\xi-1|\geq
|\zeta-1|.$ We have
\begin{eqnarray}\label{eq97}
\frac{|\psi_{\delta}(\xi)f(\xi)-\psi_{\delta}(\zeta)f(\zeta)|}{\omega(|\xi-\zeta|)}\leq
|\psi_{\delta}(\xi)|\frac{|f(\xi)-f(\zeta)|}{\omega(|\xi-\zeta|)}+|f(\zeta)|\frac{|\psi_{\delta}(\xi)-\psi_{\delta}(\zeta)|}{\omega(|\xi-\zeta|)}.
\end{eqnarray}
Suppose first that $|\xi-\zeta|\geq |\zeta-1|.$ Then
\begin{eqnarray}\nonumber
|f(\zeta)|\frac{|\psi_{\delta}(\xi)-\psi_{\delta}(\zeta)|}{\omega(|\xi-\zeta|)}&\leq&
2\frac{|f(\zeta)|}{\omega(|\zeta-1|)}\\\label{eq98}&=&o(1)\qquad(\mbox{as}\ |\xi-\zeta|\longrightarrow 0).
\end{eqnarray}
Now if $|\xi-\zeta|\leq |\zeta-1|,$ then $|z-1|\geq |\zeta-1|$ for every $z\in[\xi,\zeta].$ We obtain
\begin{eqnarray}\nonumber
|f(\zeta)|\frac{|\psi_{\delta}(\xi)-\psi_{\delta}(\zeta)|}{\omega(|\xi-\zeta|)}&=&
|f(\zeta)|\frac{|\psi_{\delta}(\xi)-\psi_{\delta}(\zeta)|}{|\xi-\zeta|}\frac{|\xi-\zeta|}{\omega(|\xi-\zeta|)}\\\nonumber&=& |f(\zeta)||\psi^{'}_{\delta}(z)|\frac{|\xi-\zeta|}{\omega(|\xi-\zeta|)}\qquad(z\in[\xi,\zeta])\\\nonumber&\leq& \frac{|f(\zeta)|}{\omega(|\zeta-1|)}\frac{\omega(|\zeta-1|)}{|\zeta-1|}\frac{|\xi-\zeta|}{\omega(|\xi-\zeta|)}\\\label{eq99}&=&o(1)\qquad(\mbox{as}\ |\xi-\zeta|\longrightarrow 0).
\end{eqnarray}
From \eqref{eq97}, \eqref{eq98} and \eqref{eq99} we deduce
$$\frac{|\psi_{\delta}(\xi)f(\xi)-\psi_{\delta}(\zeta)f(\zeta)|}{\omega(|\xi-\zeta|)}
=o(1)\qquad(\mbox{as}\ |\xi-\zeta|\longrightarrow 0),$$ uniformly
with respect to $\delta>0.$ So the result follows by applying Lemma
\ref{shilbert} to the family of functions $\psi_{\delta}f,$
$\delta>0.$ This completes the proof of the lemma.
\end{proof}

We denote by $K^c$ the complement in $\T$ of the subset $K$ of $\T.$ For a closed subset $E$ of $\T,$ we have $E^c=\bigcup\limits_{n\in\N}(a_n,b_n),$ where $(a_n,b_n)\subset E^c \mbox{ and }a_n,b_n\in E.$ We define $\Omega_{_{E}}$ to be the family of all the
unions of arcs $(a_n,b_n),$ where $(a_n,b_n)\subset E^c$ and $a_n,b_n\in E.$
For a given function $f$ in the disk algebra and $\Gamma\in\Omega_{_{E}}$ ($E$ is a closed subset of $\T$), let us define the
outer function $f_{_{\Gamma}}\in\cH^{\infty}(\D)$ associated to the
outer factor of $f$   by
$$f_{_{\Gamma}}(z):=\exp\Bigm\{\frac{1}{2\pi}\int_{\Gamma}\frac{e^{i\theta}+z}{e^{i\theta}-z}\log|f(e^{i\theta})|d\theta
\Bigm\}\qquad(z\in\D).$$
Then, we assert

\begin{lem}\label{fal}Let $\cI$ be a closed ideal of $\cL,$  $f\in\cL$  an outer function and let $h\in\cJ_{{E_{_\cI}}}$ a function such that $hf\in\cI.$
Let $\Gamma\in\Omega_{_{E_{_\cI}}}$  be such that $\T\setminus\overline{\Gamma}$ is union of a finite number of arcs $(a,b)\subseteq\T\setminus E_{_\cI}$ $(a,b\in E_{_\cI})$. If
$hf_{_\Gamma}\in\cL,$ then  $hf_{_\Gamma}\in\cI.$
\end{lem}

\begin{proof} For simplicity we suppose that $\T\setminus\overline{\Gamma}=(a,b):=\gamma,$ where
$a,b\in E_{_\cI}$ and $(a,b)\subseteq\T\setminus E_{_\cI}.$
Let $\varepsilon>0$ be such that $\gamma_{\varepsilon}:=(a e^{i\varepsilon},b e^{-i\varepsilon})\subset\gamma.$
For
$\delta>0,$ we set
$$\displaystyle  \phi_{_{\delta,\varepsilon}}(z):=\Big(\frac{z\overline{a}e^{-i\varepsilon}-1}
{z\overline{a}e^{-i\varepsilon}-1-\delta}\Big)
\Big(\frac{z\overline{b}e^{i\varepsilon}-1}{z\overline{b}e^{i\varepsilon}-1-\delta}\Big)\qquad(z\in\D).$$
It is clear that $\phi_{_{\delta,\varepsilon}}\in\cL$ and that
$\phi_{_{\delta,\varepsilon}}(a
e^{i\varepsilon})=\phi_{_{\delta,\varepsilon}}(b
e^{-i\varepsilon})=0.$ Then, according to Proposition \ref{prop1}-1
(below), we see that the function $f_{\gamma_{\varepsilon}}$
multiplied by the square of $\phi_{_{\delta,\varepsilon}}$ belongs
to $\cL$, i.e.,
$\phi^{2}_{_{\delta,\varepsilon}}f_{\gamma_{\varepsilon}}\in\cL.$
Similarly, we get $\phi^{2}_{_{\delta,\varepsilon}}
f_{\gamma^{c}_{\varepsilon}}\in\cL.$ Now, for $\pi:\cL\to \cL/\cI$
being the canonical quotient map, it follows
\begin{eqnarray*} 0&=&\pi\Big(\phi^{4}_{_{\delta,\varepsilon}}hf\Big)\\
&=&\pi\Big(h\phi^{2}_{_{\delta,\varepsilon}}f_{_{\gamma^{c}_{\varepsilon}}}\Big)\pi
\Big(\phi^{2}_{_{\delta,\varepsilon}}f_{_{\gamma_{\varepsilon}}}\Big).
\end{eqnarray*}
Since, the function
$\phi^{2}_{_{\delta,\varepsilon}}f_{_{\gamma_{\varepsilon}}}$ is
invertible in the quotient algebra $\cL/\cI,$ then
$$h\phi^{2}_{_{\delta,\varepsilon}}f_{_{\gamma^{c}_{\varepsilon}}}\in\cI\qquad(\delta,\varepsilon>0).$$
Using the fact
$$\lim\limits_{\varepsilon\longrightarrow0}\|\phi^{2}_{_{\delta,\varepsilon}}
f_{_{\gamma^{c}_{\varepsilon}}}-\phi^{2}_{_{\delta,0}}f_{_{\Gamma}}\|_{\omega}=0,$$
we can check
$$\phi^{2}_{_{\delta,0}}hf_{_{\Gamma}}\in\cI\qquad(\delta>0).$$
Now, since $hf_{_{\Gamma}}\in\cL$ and $hf_{_{\Gamma}}(a)=hf_{_{\Gamma}}(b)=0,$ we can deduce from
Lemma \ref{lem3} that
$$\lim\limits_{\delta\longrightarrow0}\|\phi^{2}_{_{\delta,0}}hf_{_{\Gamma}}-hf_{_{\Gamma}}\|_{\omega}=0.$$
So  $hf_{_{\Gamma}}\in\cI.$ This completes the proof of the lemma.
\end{proof}

Let $\cI$ be a closed ideal of the algebra $\cL.$ We have
$U_{_\cI}=B_{_\cI}S_{_\cI}.$ The inner function $B_{_\cI}$ is the
usual Blashke product associated to the zero set $Z_{_\cI}\cap\D,$
where $Z_{_\cI}:=\{z\in\overline{\D}\ :\ f(z)=0 \mbox{ for all
}f\in\cI\}.$ The positive singular measure $\mu_{_\cI}$ associated
to the singular inner function $S_{_\cI}$ is the greatest common
divisor of all $\mu_{_f},$ $f\in\cI.$  Note that
$\mathrm{supp}(\mu_{_\cI})$ is included in $E_{_\cI}.$ For a subset $K$ of $T,$  we set $$\big(S_{_f}\big)_{K}(z):=\exp\Big\{-\frac{1}{2\pi}\int_{K}\frac{e^{i\theta}+z}{e^{i\theta}-z}d\mu_{_f}(\theta)\Big\}\qquad(f\in\cA(\D)).$$

\begin{lem}\label{lem2} Let $f$ be a function in a closed ideal $\cI$ of the algebra $\cL.$
Then $B_{_\cI}\big(S_{_f}\big)_{E_{_\cI}}O_{_f}$ belongs to $\cI.$
\end{lem}

\begin{proof}
Let $f\in\cI.$ Define $B_{f,n}$ and $B_{\cI,n}$ to be respectively the Blashke product with zeros $Z_{_f}\cap\D_{n}$ and  $Z_{_\cI}\cap\D_{n},$ where $\D_{n}:=\{z\in\D\text{ : }|z|<\frac{n-1}{n},\ n\in\N\}.$ Fix $n\in\N.$ The function $B_{f,n}/B_{\cI,n}$ is invertible in the quotient algebra $\cL/\cJ_n,$ where $\cJ_n:=\{g\in\cL\ :\ gB_{\cI,n}\in\cI\}.$
Then $f/B_{f,n}\in\cJ_n.$ It follows that $B_{\cI,n}(f/B_{f,n})\in\cI.$
It is clear that
$$\lim\limits_{n\longrightarrow+\infty}\|B_{\cI,n}(f/B_{f,n})-B_{_\cI}S_{_f}O_{_f}\|_{\infty}=0.$$
By using Corollary \ref{fprc} in appendix B (F-property of $\cL$), we obtain
$$\lim\limits_{n\longrightarrow+\infty}\|B_{\cI,n}(f/B_{f,n})-B_{_\cI}S_{_f}O_{_f}\|_{\omega}=0.$$
So, $B_{_\cI}S_{_f}O_{_f}\in\cI.$ Let $\varepsilon>0$ be such that $\gamma_{\varepsilon}:=(a e^{i\varepsilon},b e^{-i\varepsilon})\subset\gamma:=(a,b),$ where
$a,b\in E_{_\cI}$ and $(a,b)\subseteq\T\setminus E_{_\cI}.$ By using the F-property of $\cL,$ the function $B_{_\cI}\big(S_{_f}\big)_{\gamma^{c}_{\varepsilon}}O_{_f}$
belongs to $\cL.$ We set
$$\displaystyle \phi_{_{\varepsilon}}(z):=(z\overline{a}e^{-i\varepsilon}-1)
(z\overline{b}e^{i\varepsilon}-1)\qquad(z\in\D).$$
From Proposition \ref{prop1}-1 below, the function $g_{_\varepsilon}:=\phi^{2}_{_{\varepsilon}}f_{\gamma_{\varepsilon}}$ belongs to $\cL.$ Since $\mathrm{supp}(\mu_{(S_{_f})_{_{\gamma_{\varepsilon}}}})\subseteq E_{{g_{_\varepsilon}}},$ then we deduce from Lemma \ref{lem1} that $\big(S_{_f}\big)_{\gamma_{\varepsilon}}g^{2}_{_\varepsilon}\in\cL.$  We have
\begin{eqnarray*}
0=\pi\Big(g^{2}_{_\varepsilon}B_{_\cI}S_{_f}O_{_f}\Big)=
\pi\Big(\big(S_{_f}\big)_{\gamma_{\varepsilon}}g^{2}_{_\varepsilon}\Big)
\times \pi\Big(B_{_\cI}\big(S_{_f}\big)_{\gamma^{c}_{\varepsilon}}O_{_f}\Big)
\end{eqnarray*}
where $\pi:\cL\to \cL/\cI,$ is the canonical quotient map.
The function $\big(S_{_f}\big)_{\gamma_{\varepsilon}}g^{2}_{_\varepsilon}$ is invertible in the
quotient algebra $\cL/\cI,$ then
$B_{_\cI}\big(S_{_f}\big)_{\gamma^{c}_{\varepsilon}}O_{_f}\in\cI.$ It is clear that
$$\lim\limits_{\varepsilon\longrightarrow 0}\|B_{_\cI}\big(S_{_f}\big)_{\gamma^{c}_{\varepsilon}}O_{_f}-B_{_\cI}\big(S_{_f}\big)_{\gamma^{c}}O_{_f}\|_{\infty}=0.$$
Then, using Corollary \ref{fprc} in appendix B, we obtain $$\lim\limits_{\varepsilon\longrightarrow 0}\|B_{_\cI}\big(S_{_f}\big)_{\gamma^{c}_{\varepsilon}}O_{_f}-B_{_\cI}\big(S_{_f}\big)_{\gamma^{c}}O_{_f}\|_{\omega}=0.$$
So $B_{_\cI}\big(S_{_f}\big)_{\gamma^{c}}O_{_f}\in\cI.$ Similarly we can prove that
$B_{_\cI}\big(S_{_f}\big)_{\Gamma_{_N}^{c}}O_{_f}\in\cI,$ where $\Gamma_{_N}:=\bigcup\limits_{n\leq N}(a_n,b_n)\in\Omega_{_{E_{_\cI}}}.$
We have
$$\lim\limits_{N\longrightarrow +\infty}\|B_{_\cI}\big(S_{_f}\big)_{\Gamma_{_N}^{c}}O_{_f}-B_{_\cI}\big(S_{_f}\big)_{E_{_\cI}}O_{_f}\|_{\infty}=0.$$
Using again Corollary \ref{fprc} we deduce
$$\lim\limits_{N\longrightarrow +\infty}\|B_{_\cI}\big(S_{_f}\big)_{\Gamma_{_N}^{c}}O_{_f}-B_{_\cI}\big(S_{_f}\big)_{E_{_\cI}}O_{_f}\|_{\omega}=0.$$
Then $B_{_\cI}\big(S_{_f}\big)_{E_{_\cI}}O_{_f}\in\cI.$ This proves the lemma.
\end{proof}

\subsection{\bf Proof of Theorem \ref{g2byf} }

For the proof of Theorem \ref{g2byf}, we need the following
proposition about approximation of functions in $\cL$

\begin{prop}\label{prop1}Let $\omega$ be a modulus of continuity satisfying the condition \eqref{2condition}. Let $f\in\cL$ be a function such that $\|f\|_{\omega}\leq 1$ and  $E$  a closed subset of $\T.$
Let $g\in\cJ_{{E}}$ be an outer function and $S$ singular inner
function such that $\mathrm{supp}(\mu_{_S})\subseteq E_{_g}.$ Then
\begin{itemize}
\item[$1.$] The functions $Sg^2$ and $Sg^2 f_{_{\Gamma^{c}_{_N}}}$ belong to $\cL,$ for every $N\in\N,$
\item[$2.$] We have $\lim\limits_{N\longrightarrow +\infty}\|Sg^2 f_{_{\Gamma^{c}_{_N}}}-Sg^2\|_{\omega}=0,$
\end{itemize}
where $\Gamma_{_N}:=\bigcup\limits_{n\leq N}(a_n,b_n)\in\Omega_{_{E}}.$ If  moreover $\omega$  satisfy the  stronger condition \eqref{3condition}, then
\begin{itemize}
\item[$1'.$] The functions $Sg$ and $Sg f_{_{\Gamma^{c}_{_N}}}$ belong to $\cL,$ for every $N\in\N,$
\item[$2'.$] We have $\lim\limits_{N\longrightarrow +\infty}\|Sg f_{_{\Gamma^{c}_{_N}}}-Sg\|_{\omega}=0.$
\end{itemize}
\end{prop}

\begin{proof} From Lemma \ref{lem1} we have $Sg^2\in\cL.$ It is clear that
$\{Sg^2 f_{_{\Gamma^{c}_{_N}}}\}_{N\in\N}$ is a sequence of functions in the disk algebra. Let us note that if
\begin{equation}\label{if}
\frac{|S(z)g^2(z)f_{_{\Gamma}}(z)-S(w)g^2(w)f_{_{\Gamma}}(w)|}{\omega(|z-w|)}=
o(1) \qquad(\text{as } |z-w|\rightarrow 0 ),
\end{equation}
uniformly with respect to $\Gamma\in\Omega_{_{E}},$ then  assertion
1 follows immediately. Furthermore,
  assertion 2  can be deduced by applying  Lemma \ref{shilbert}.\\
  Thus, it suffices to show only \eqref{if}. For this, we fix $\Gamma\in\Omega_{_{E}}$ and we let $\xi$ and $\zeta$ be two distinct points in $\T$ such that $d(\xi,E)\geq d(\zeta,E).$  It is clear that
\begin{eqnarray*}
&&\frac{|S(\xi)g^{2}(\xi)f_{_\Gamma}(\xi)-S(\zeta)g^{2}(\zeta)f_{_\Gamma}(\zeta)|}{\omega(|\xi-\zeta|)}
\\&\leq&|f_{_\Gamma}(\xi)|\frac{|S(\xi)g^{2}(\xi)-S(\zeta)g^{2}(\zeta)|}{\omega(|\xi-\zeta|)}+|g^{2}(\zeta)|\frac{|f_{_\Gamma}(\xi)-f_{_\Gamma}(\zeta)|}{\omega(|\xi-\zeta|)}.
\end{eqnarray*}
Since $Sg^{2}\in\cL$ (by Lemma \ref{lem1}), then the proof of \eqref{if}  reduces to
\begin{equation}\label{gbyfp}
|g^{2}(\zeta)|\frac{|f_{_\Gamma}(\xi)-f_{_\Gamma}(\zeta)|}{\omega(|\xi-\zeta|)}= o(1) \qquad(\text{as } |\xi-\zeta|\longrightarrow 0 ).
\end{equation}

 \textit{Case} 1. For
$\displaystyle|\xi-\zeta|\geq \big(\frac{d(\zeta,E)}{2}\big)^{2},$ we have
\begin{eqnarray}\nonumber
|g^{2}(\zeta)|\frac{|f_{_\Gamma}(\xi)-f_{_\Gamma}(\zeta)|}{\omega(|\xi-\zeta|)}&\leq&
2\frac{|g^{2}(\zeta)|}{\omega(|\xi-\zeta|)}\\\nonumber
&\leq& 8\eta_{2}^{-1}
\Big(\frac{|g(\zeta)|}{\omega(d(\zeta,E))}\Big)^{2}
\\&=&\label{gbyf1}o(1) \qquad(\text{as } |\xi-\zeta|\longrightarrow 0 ).
\end{eqnarray}

\textit{Case} 2. For $\displaystyle|\xi-\zeta|\leq \big(\frac{d(\zeta,E)}{2}\big)^{2}$ with
$\zeta\notin\Gamma.$ It follows $[\xi,\zeta]\subset \big(\Gamma\cup E\big)^{c}.$ Then $z\notin\Gamma$ and $ d(z,E)\geq d(\zeta,E)$ for
every $z\in[\xi,\zeta].$ There is $z\in[\xi,\zeta]$ such that
$$ \frac{|f_{_\Gamma}(\xi)-f_{_\Gamma}(\zeta)|}{|\xi-\zeta|}=|f_{_\Gamma}^{'}(z)|\leq \frac{1}{\pi}\int_\Gamma\frac{|\log|f(e^{i\theta})||}{|e^{i\theta}-z|^2}d\theta\leq \frac{ c_{_f}}{d^2(z,E)}.$$
It follows  $$ \frac{|f_{_\Gamma}(\xi)-f_{_\Gamma}(\zeta)|}{|\xi-\zeta|}\leq \frac{c_{_f}}{d^2(\zeta,E)}.$$ Therefore, we have
\begin{eqnarray}\nonumber
|g^{2}(\zeta)|\frac{|f_{_\Gamma}(\xi)-f_{_\Gamma}(\zeta)|}{\omega(|\xi-\zeta|)}&=&
\frac{|g^{2}(\zeta)||\xi-\zeta|}{\omega(|\xi-\zeta|)}\frac{|f_{_\Gamma}(\xi)-f_{_\Gamma}(\zeta)|}{|\xi-\zeta|}\\\nonumber&\leq&
c_f\frac{|g^{2}(\zeta)||\xi-\zeta|}{\omega(|\xi-\zeta|)d^2(\zeta,E)}\\\nonumber&\leq&
c_f\Big(\frac{|g(\zeta)|}{\omega(d(\zeta,E))}\Big)^2 \frac{\omega(d^2(\zeta,E))|\xi-\zeta|}{\omega(|\xi-\zeta|)d^2(\zeta,E)}\\\label{eee}&=&
o(1) \qquad(\text{as } |\xi-\zeta|\longrightarrow 0 ).
\end{eqnarray}

\textit{Case} 3. In this case let us assume that $\displaystyle|\xi-\zeta|\leq
\big(\frac{d(\zeta,E)}{2}\big)^{2}$ and $\zeta\in\Gamma$.  It follows that $\xi\in\Gamma,$ and therefore $|f^{-1}_{_{\T\setminus\Gamma}}(\xi)|=|f^{-1}_{_{\T\setminus\Gamma}}(\zeta)|=1.$ From
\begin{eqnarray}\nonumber
&&f_{_\Gamma}(\xi)-f_{_\Gamma}(\zeta)\\\nonumber&=&
f(\xi)f^{-1}_{_{\T\setminus\Gamma}}(\xi)-f(\zeta)f^{-1}_{_{\T\setminus\Gamma}}(\zeta)\\\nonumber&=&
f^{-1}_{_{\T\setminus\Gamma}}(\xi)\big(f(\xi)-f(\zeta)\big)+f(\zeta)\big(f^{-1}_{_{\T\setminus\Gamma}}(\xi)-f^{-1}_{_{\T\setminus\Gamma}}(\zeta)\big)\\\nonumber&=&
f^{-1}_{_{\T\setminus\Gamma}}(\xi)\big(f(\xi)-f(\zeta)\big)\\\label{revien}
&&-f(\zeta)f^{-1}_{_{\T\setminus\Gamma}}(\xi)f^{-1}_{_{\T\setminus\Gamma}}(\zeta)
\big(f_{_{\T\setminus\Gamma}}(\xi)-f_{_{\T\setminus\Gamma}}(\zeta)\big)
\end{eqnarray}
combined with \eqref{eee},  we deduce
\begin{eqnarray}\nonumber&&
|g^{2}(\zeta)|\frac{|f_{_\Gamma}(\xi)-f_{_\Gamma}(\zeta)|}{\omega(|\xi-\zeta|)}\\\nonumber&\leq&
|g^{2}(\zeta)|\frac{|f(\xi)-f(\zeta)|}{\omega(|\xi-\zeta|)}+|f(\zeta)||g^{2}(\zeta)|\frac{|f_{_{\T\setminus\Gamma}}(\xi)-f_{_{\T\setminus\Gamma}}(\zeta)|}{\omega(|\xi-\zeta|)}\\\label{gbyfd}&=&
o(1)\qquad(\text{as } |\xi-\zeta|\longrightarrow 0 ).
\end{eqnarray}
From inequalities \eqref{gbyf1}, \eqref{eee} and \eqref{gbyfd}, we
see that \eqref{gbyfp} holds. If moreover $\omega$ satisfy the condition \eqref{3condition}, we can deduce similarly the assertions 1' and 2'. This completes the proof of the
proposition.
\end{proof}

\textit{\bf Proof of Theorem \ref{g2byf}:} Let $\cI$ be a closed
ideal of the algebra $\cL$ such that $U_{_{\cI}}\equiv1$. Let $g$ be
an outer function in $\cJ_{E_{_{\cI}}}$ and consider $f\in\cI.$
From Lemma \ref{lem2}, we deduce that
$\big(S_{_f}\big)_{E_{_{\cI}}}O_{_f}\in\cI.$ Then,
$(S_{_f})_{E_{_{\cI}}}g^2 O_{_f}\in\cI.$  From Lemma \ref{fal} and Proposition \ref{prop1}-1, we
deduce that $\big(S_{_f}\big)_{E_{_{\cI}}}g^2f_{_{\Gamma^{c}_{_N}}}\in\cI,$
for every $N\in\N,$  where $\Gamma_{_N}:=\bigcup\limits_{n\leq N}(a_n,b_n)\in\Omega_{_{E_{_{\cI}}}}.$ Thus, by  Proposition \ref{prop1}-2,
$\big(S_{_f}\big)_{E_{_{\cI}}}g^2\in\cI.$ Now choose a sequence of
functions  $\{f_n\}_{n\in\N}\subset\cI$ such that the greatest
common divisor of the inner parts of $f_n$ is equal to 1 and such that $(1/2\pi)\int_{0}^{2\pi}d\mu_{_{f_n}}(\theta)\leq1.$ This
infers that $k_n:=\big(S_{_{f_n}}\big)_{E_{_{\cI}}}g^2\in\cI,$ for every
$n\in\N.$ By Lemma \ref{lem1} we have
\begin{equation*}
\frac{|k_n(\xi)-k_n(\zeta)|}{\omega(|\xi-\zeta|)}=o(1)\qquad (\mbox{as } |\xi-\zeta|\longrightarrow 0),
\end{equation*}
uniformly with respect to $n.$ Using the fact
$$\lim\limits_{n\longrightarrow+\infty}\|k_n-g^{2}\|_{\infty}=0$$
and Lemma \ref{shilbert}, we deduce
$$\lim\limits_{n\longrightarrow+\infty}\|k_n-g^{2}\|_{\omega}=0.$$
So $g^2\in\cI.$ This finishes the proof of the theorem.

\subsection{\bf  Proof of Theorem \ref{gbyg2}}

We begin by proving the following proposition

\begin{prop}\label{prop3}Let $\omega$ be a modulus of continuity satisfying the condition \eqref{cond}.
Fix $1<\rho\leq2$. Let $g\in\cL$ be a function such that
$\|g\|_{\omega}\leq 1$  and let
 $E$ be a closed subset of $E_{_g}.$  Then
\begin{itemize}
  \item[$1.$] The functions $U_{_g}O^\rho_g$ and $U_{g}O^\rho_gg_{_{\Gamma^{c}_{_N}}}$ belong to $\cL,$ for every $N\in\N,$
  \item[$2.$] We have $ \lim\limits_{N\longrightarrow +\infty}\|U_{_g}O^\rho_gg_{_{\Gamma^{c}_{_N}}}-U_{_g}O^\rho_g\|_{\omega}=0,$
\end{itemize}
 where $\Gamma_{_N}:=\bigcup\limits_{n\leq N}(a_n,b_n)\in\Omega_{_{E}}.$
\end{prop}

\begin{proof}By using the F-property of $\cL$ and the following
\begin{eqnarray*}
&&U_{_g}(\xi)O^\rho_g(\xi)-U_{_g}(\zeta)O^\rho_g(\zeta)\\&=&
O^{\rho-1}_g(\xi)\big(U_{_g}(\xi)O_g(\xi)-U_{_g}(\zeta)O_g(\zeta)\big)-
U_{_g}(\zeta)O^{\rho-1}_g(\xi)\big(O_g(\xi)-O_g(\zeta)\big)\\&+&
U_{_g}(\zeta)\big(O^{\rho}_g(\xi)-O^{\rho}_g(\zeta)\big)\qquad(\xi,\zeta\in\T),
\end{eqnarray*}
the function $U_{_g}O_g^{\rho}$ belongs to $\cL.$ It is clear that
$\{U_{g}O^\rho_gg_{_{\Gamma^{c}_{_N}}}\}_{N\in\N}$ is sequence of functions in the disk algebra and
$$\lim\limits_{N\longrightarrow +\infty}\|U_{g}O^\rho_gg_{_{\Gamma^{c}_{_N}}}-U_{g}O^\rho_g\|_{\infty}=0.$$
Note that if
\begin{equation}\label{iff}
\frac{|U_{g}O^\rho_gg_{_\Gamma}(\xi)-U_{g}O^\rho_gg_{_\Gamma}(\zeta)|}{\omega(|\xi-\zeta|)}=
o(1) \qquad(\text{as } |\xi-\zeta|\rightarrow 0 ),
\end{equation}
uniformly with respect to $\Gamma\in\Omega_{_{E}},$ then  assertion
1 holds as well as assertion 2, which follows from  \eqref{iff}
combined with Lemma \ref{shilbert}.\\ Below, we have to prove
\eqref{iff}. Let $\xi,\zeta\in\T$ be two different points such that $d(\xi,E)\geq d(\zeta,E).$  It is clear that
\begin{eqnarray*}
&&\frac{|U_{_g}(\xi)O^\rho_g(\xi)g_{_\Gamma}(\xi)-U_{_g}(\zeta)O^\rho_g(\zeta)g_{_\Gamma}(\zeta)|}{\omega(|\xi-\zeta|)}
\\&\leq&|g_{_\Gamma}(\xi)|\frac{|U_{_g}(\xi)O^\rho_g(\xi)-U_{_g}(\zeta)O^\rho_g(\zeta)|}{\omega(|\xi-\zeta|)}+|g(\zeta)|^\rho\frac{|g_{_\Gamma}(\xi)-g_{_\Gamma}(\zeta)|}{\omega(|\xi-\zeta|)}.
\end{eqnarray*}
Then, to prove \eqref{iff} it is sufficient to prove that
\begin{equation}\label{gbygp}
|g(\zeta)|^\rho\frac{|g_{_\Gamma}(\xi)-g_{_\Gamma}(\zeta)|}{\omega(|\xi-\zeta|)}=o(1) \qquad(\text{as } |\xi-\zeta|\longrightarrow 0 ).
\end{equation}
\begin{itemize}
  \item [$1.$] We suppose that $\displaystyle|\xi-\zeta|\geq \big(\frac{d(\zeta,E)}{2}\big)^\rho,$ we obtain
\begin{eqnarray}\nonumber
|g(\zeta)|^\rho\frac{|g_{_\Gamma}(\xi)-g_{_\Gamma}(\zeta)|}{\omega(|\xi-\zeta|)}&\leq&
2\frac{|g(\zeta)|^\rho}{\omega(|\xi-\zeta|)}\\\nonumber
&\leq& 2^{\rho+1}\eta^{-1}_{\rho}
\Big(\frac{|g(\zeta)|}{\omega(d(\zeta,E))}\Big)^\rho
\\\label{aaa}&=&o(1) \qquad(\text{as } |\xi-\zeta|\longrightarrow 0 ).
\end{eqnarray}
\item [$2.$] In this case we suppose that $\displaystyle|\xi-\zeta|\leq\big(\frac{d(\zeta,E)}{2}\big)^\rho$
and that  $\zeta\notin\Gamma.$ It follows $[\xi,\zeta]\subset\T\setminus E.$ Then $z\notin\Gamma$ and $ |z-e^{i\theta}|\geq \frac{1}{2}|\zeta-e^{i\theta}|$ for
every $z\in[\xi,\zeta]$ and for every $e^{i\theta}\in\Gamma.$ There is $z\in[\xi,\zeta]$ such that
$$\frac{|g_{_\Gamma}(\xi)-g_{_\Gamma}(\zeta)|}{|\xi-\zeta|}=|g^{'}_{_\Gamma}(z)|\leq a_{_\Gamma}(z),$$
where $\displaystyle
a_{_\Gamma}(z):=\frac{1}{\pi}\int_\Gamma\frac{|\log|g(e^{i\theta})||}{|e^{i\theta}-z|^2}d\theta.$
Since  $a_{_\Gamma}(z)\leq 4\ a_{_\Gamma}(\zeta),$ then
$$\frac{|g_{_\Gamma}(\xi)-g_{_\Gamma}(\zeta)|}{|\xi-\zeta|}\leq 4\ a_{_\Gamma}(\zeta).$$
\begin{itemize}
\item [$2.1.$] First we suppose that $\displaystyle a_{_\Gamma}(\zeta)\leq
\frac{1}{d^\rho(\zeta,E)}.$ It is clear that
\begin{eqnarray*}
\frac{|g(\zeta)|^{\rho}|\xi-\zeta|}{d^\rho(\zeta,E)\omega(|\xi-\zeta|)}&=&
\Big(\frac{|g(\zeta)|}{\omega(d(\zeta,E))}\Big)^{\rho}\Big(\frac{\omega(d(\zeta,E))}{d(\zeta,E)}\Big)^\rho
\frac{|\xi-\zeta|}{\omega(|\xi-\zeta|)}\\&\leq&
\inf\Big\{\eta^{-1}_{\rho}\Big(\frac{|g(\zeta)|}{\omega(d(\zeta,E))}\Big)^{\rho},
\Big(\frac{\omega(d(\zeta,E))}{d(\zeta,E)}\Big)^{\rho}\frac{|\xi-\zeta|}{\omega(|\xi-\zeta|)}\Big\}\\
&=&o(1)\qquad (\mbox{as } |\xi-\zeta|\longrightarrow 0).
\end{eqnarray*}
Therefore we obtain,
\begin{eqnarray}\nonumber
|g(\zeta)|^{\rho}\frac{|g_{_\Gamma}(\xi)-g_{_\Gamma}(\zeta)|}{\omega(|\xi-\zeta|)}&\leq&
\frac{|g(\zeta)|^{\rho}|\xi-\zeta|}{\omega(|\xi-\zeta|)}\frac{|g_{_\Gamma}(\xi)-g_{_\Gamma}(\zeta)|}{|\xi-\zeta|}
\\\nonumber&\leq&
4\frac{|g(\zeta)|^{\rho}|\xi-\zeta|}{d^\rho(\zeta,E)\omega(|\xi-\zeta|)}\\
&=&\label{1est}o(1) \qquad(\text{as } |\xi-\zeta|\longrightarrow 0).
\end{eqnarray}
\item [$2.2.$] Next we suppose that $\displaystyle a_{_\Gamma}(\zeta)\geq
\frac{1}{d^{\rho-1}(\zeta,E)|\xi-\zeta|^{1/\rho}}.$ Set
$\lambda_\zeta:=1-|\xi-\zeta|^{1/\rho}.$ Then
\begin{eqnarray*}
|g(\lambda_\zeta\zeta)|&=&\exp\Bigm\{\frac{1}{2\pi}\int_{0}^{2\pi}\frac{1-\lambda_\zeta^{2}
}{|e^{i\theta}- \lambda_\zeta\zeta|^{2}}\log|g(e^{i\theta})|d\theta\Bigm\}\\
&\leq&
\exp\Bigm\{\frac{1}{2\pi}\int_{\Gamma}\frac{|\xi-\zeta|^{1/\rho}}{|e^{i\theta}-\lambda_\zeta\zeta|^{2}}\log|g(e^{i\theta})|d\theta\Bigm\}
\\  &\leq&\exp\Bigm\{-\frac{1}{4}|\xi-\zeta|^{1/\rho}a_{_\Gamma}(\zeta)\Bigm\}
\\ &\leq& \exp\Bigm\{-\frac{1}{4d^{\rho-1}(\zeta,E)}\Bigm\}.
\end{eqnarray*}
It is clear that
$\displaystyle|g(\zeta)|^{\rho}\leq 2^{\rho-1}|g(\zeta)-g(\lambda_\zeta\zeta)|^{\rho}+2^{\rho-1}|g(\lambda_\zeta\zeta)|^{\rho}.$
We obtain
\begin{eqnarray}&&\nonumber
|g(\zeta)|^{\rho}\frac{|g_{_\Gamma}(\xi)-g_{_\Gamma}(\zeta)|}{\omega(|\xi-\zeta|)}\\\nonumber&\leq&
2^{\rho-1}\frac{|g(\zeta)-g(\lambda_\zeta\zeta)|^\rho}{\omega(|\xi-\zeta|)}|g_{_\Gamma}(\xi)-g_{_\Gamma}(\zeta)|\\\nonumber&&
+2^{\rho-1}|g(\lambda_\zeta\zeta)|^{\rho}\frac{|\xi-\zeta|}{\omega(|\xi-\zeta|)}\frac{|g_{_\Gamma}(\xi)-g_{_\Gamma}(\zeta)|}{|\xi-\zeta|}\\\nonumber&\leq&
2^\rho\eta^{-1}_{\rho}o(1)\\\nonumber&&
+2^{\rho-1}\frac{|\xi-\zeta|}{\omega(|\xi-\zeta|)}\exp\Bigm\{-\frac{\rho}{4d^{\rho-1}(\zeta,E)}\Bigm\}\frac{\int_{\T}|\log|g(e^{i\theta})||d\theta}{d^2(\zeta,E)}\\\label{2est}&=&
o(1)\qquad(\mbox{as }|\xi-\zeta|\longrightarrow 0).
\end{eqnarray}
\item [$2.3.$] Now we suppose that $\displaystyle \frac{1}{d^{\rho}(\zeta,E)}\leq a_{_\Gamma}(\zeta)\leq
\frac{1}{d^{\rho-1}(\zeta,E)|\xi-\zeta|^{1/\rho}}.$ Set
$\displaystyle \mu_\zeta:=1-\frac{1}{a_{_\Gamma}(\zeta)d^{\rho-1}(\zeta,E)}.$ Then
$|\xi-\zeta|^{1/\rho}\leq 1-\mu_\zeta\leq d(\zeta,E).$ It follows that
\begin{eqnarray*}\frac{\omega(1-\mu_\zeta)}{1-\mu_\zeta}\frac{|\xi-\zeta|}{\omega(|\xi-\zeta|)}
\Big(\frac{\omega(d(\zeta,E))}{d(\zeta,E)}\Big)^{\rho-1}&\leq&
\Big(\frac{\omega(1-\mu_\zeta)}{1-\mu_\zeta}\Big)^{\rho}\frac{|\xi-\zeta|}{\omega(|\xi-\zeta|)}\\
&\leq& \eta^{-1}_{\rho},
\end{eqnarray*}
and  $$\frac{\omega(1-\mu_\zeta)}{1-\mu_\zeta}\frac{|\xi-\zeta|}{\omega(|\xi-\zeta|)}\leq
\eta^{-(1/\rho)}_{\rho}\Big(\frac{|\xi-\zeta|}{\omega(|\xi-\zeta|)}\Big)^{\frac{\rho-1}{\rho}}.$$
 Then
\begin{eqnarray*}
&&\frac{\omega(1-\mu_\zeta)}{1-\mu_\zeta}\frac{|\xi-\zeta|}{\omega(|\xi-\zeta|)}
\Big(\frac{|g(\zeta)|}{d(\zeta,E)}\Big)^{\rho-1}\\&\leq& \inf\big\{\eta^{-1}_{\rho}\Big(\frac{|g(\zeta)|}{\omega(d(\zeta,E))}\Big)^{\rho-1},
\eta^{-(1/\rho)}_{\rho}\Big(\frac{|\xi-\zeta|}{\omega(|\xi-\zeta|)}\Big)^{\frac{\rho-1}{\rho}}
\Big(\frac{|g(\zeta)|}{d(\zeta,E)}\Big)^{\rho-1}\big\}\\
&=&o(1)\quad(\mbox{as }|\xi-\zeta|\longrightarrow 0).
\end{eqnarray*}
Also we have
$\displaystyle|g(\mu_\zeta\zeta)|\leq\exp\Bigm\{-\frac{1}{4d^{\rho-1}(\zeta,E)}\Bigm\}.$
We obtain
\begin{eqnarray}&&\nonumber
|g(\zeta)|^{\rho}\frac{|g_{_\Gamma}(\xi)-g_{_\Gamma}(\zeta)|}{\omega(|\xi-\zeta|)}\\\nonumber&\leq&
|g(\zeta)|^{\rho}\frac{|\xi-\zeta|}{\omega(|\xi-\zeta|)}a_{_\Gamma}(z)\qquad(z\in[\xi,\zeta])\\\nonumber&\leq&
|g(\zeta)|^{\rho-1}\frac{|g(\zeta)-g(\mu_\zeta\zeta)||\xi-\zeta|}{(1-\mu_\zeta)\omega(|\xi-\zeta|)}(1-\mu_\zeta)a_{_\Gamma}(z)\\\nonumber&&
+
|g(\zeta)|^{\rho-1}|g(\mu_\zeta\zeta)|\frac{|\xi-\zeta|}{\omega(|\xi-\zeta|)}a_{_\Gamma}(z)\\\nonumber&\leq&
\frac{\omega(1-\mu_\zeta)}{1-\mu_\zeta}\frac{|\xi-\zeta|}{\omega(|\xi-\zeta|)}\Big(\frac{|g(\zeta)|}{d(\zeta,E)}\Big)^{\rho-1}\\\nonumber&&+
\frac{|\xi-\zeta|}{\omega(|\xi-\zeta|)}\exp\Bigm\{-\frac{1}{4d^{\rho-1}(\zeta,E)}\Bigm\}\frac{\int_{\T}|\log|g(e^{i\theta})||d\theta}{d^2(\zeta,E)}
\\\label{3est}&=& o(1)\qquad(\mbox{as }|\xi-\zeta|\longrightarrow 0).
\end{eqnarray}
\end{itemize}
Consequently,  from \eqref{1est}, \eqref{2est} and \eqref{3est}
we obtain that if $|\xi-\zeta|\leq \big(\frac{d(\zeta,E)}{2}\big)^\rho$ and $\zeta\notin\Gamma,$ then
\begin{equation}\label{4est}
|g^{\rho}(\zeta)|\frac{|g_{_\Gamma}(\xi)-g_{_\Gamma}(\zeta)|}{\omega(|\xi-\zeta|)}=o(1)\qquad(\mbox{as }|\xi-\zeta|\longrightarrow 0).
\end{equation}
\item [3.] In this case we suppose that $\displaystyle|\xi-\zeta|\leq \big(\frac{d(\zeta,E)}{2}\big)^\rho$
and that  $\zeta\in\Gamma.$ We use the equality \eqref{revien} to
transfer this case into the case $\zeta\notin\Gamma.$ Hence we use
\eqref{4est} to obtain also in this case that
\begin{equation}\label{5est}
|g^{\rho}(\zeta)|\frac{|g_{_\Gamma}(\xi)-g_{_\Gamma}(\zeta)|}{\omega(|\xi-\zeta|)}=o (1)\qquad(\mbox{as }|\xi-\zeta|\longrightarrow 0).
\end{equation}
\end{itemize}
Now \eqref{gbygp} follows from inequalities \eqref{aaa}, \eqref{4est} and \eqref{5est}.
This completes the proof of the proposition.
\end{proof}

\textit{\bf Proof of Theorem \ref{gbyg2} :} Let $\cI$ be a closed
ideal in $\cL$ and  $g$  a function in $\cL.$ Then $O_{_g}\in\cL$
and $O^{\rho}_{_g}\in\cL,$ for every $\rho>1.$ Suppose that
$U_{_g}O^2_{_g}\in\cI.$ Then $U_{_g}O^{\rho+2}_{_g}\in\cI,$ for
every $\rho>1.$  From Proposition
\ref{prop3}-1 we have $U_{_g}O^{\rho}_{_g}g_{_{\Gamma^{c}_{_N}}}\in\cL,$ and hence, $U_{_g}O^{\rho+1}_{_g}g_{_{\Gamma^{c}_{_N}}}\in\cL.$
It follows, by using Lemma \ref{fal}, that
$U_{_g}O^{\rho+1}_{_g}g_{_{\Gamma^{c}_{_N}}}\in\cI,$ for every $\rho>1$
and for every $N\in\N.$ Hence, from Proposition \ref{prop3}-2,
$U_{_g}O^{\rho+1}_{_g}\in\cI$ for every $\rho>1.$ Again we use Lemma
\ref{fal} and Proposition \ref{prop3} to deduce that
$U_{_g}O^{\rho}_{_g}\in\cI$ for every $\rho>1.$ This infers that
$g\in\cI.$

\section{\bf Appendix A. An equivalent norm in $\cL.$}

In this section we give a simple proof of the following Tamrazov's Theorem \cite{Tam}.

\begin{thm}\cite{Tam}\label{Tam}
Let $\omega$ be any arbitrary modulus of continuity and let $f$ be a function
 in $\cL(\T)\cap\cA(\D).$ Then $f$ belongs to $\cL.$
\end{thm}

We need the following key lemma.
\begin{lem}\label{key0}
Let $\omega$ be any arbitrary modulus of continuity and let $f\in\cL(\T)\cap\cA(\D)$ be a function
such that $\|f\|_{\cL(\T)}\leq 1.$  For every $\delta\geq0,$ we have
\begin{eqnarray*}&&\exp\Big\{\frac{1}{2\pi}\int_{0}^{2\pi}\frac{1-|z|^2}
{|e^{i\theta}-z|^2}\log\big(|f(e^{i\theta})-f(z/|z|)|+\delta\big)d\theta\Big\}
\\&\leq& o(\omega(1-|z|))+A\delta\qquad(\text{as }|z|\longrightarrow1),
\end{eqnarray*}
where $A$ is an absolute constant.
\end{lem}

\begin{proof} Let $0<\varepsilon<1$ and
$c_\varepsilon>0$ such that for any $\xi,\zeta\in\T$
satisfying $|\xi-\zeta|\leq c_\varepsilon,$ we have
$|f(\xi)-f(\zeta)|\leq \varepsilon\ \omega(|\xi-\zeta|) .$ Divide
$\T$ into the following three parts
\begin{eqnarray*}
\Gamma_1&:=& \big\{\xi\in\T \ :\ |\xi-z/|z||\leq 1-|z|\leq c_{\varepsilon}\big\},\\
\Gamma_2&:=& \big\{\xi\in\T \ :\ 1-|z|\leq|\xi-z/|z||\leq c_{\varepsilon}\big\},\\
\Gamma_3&:=& \big\{\xi\in\T \ :\ c_{\varepsilon}\leq|\xi-z/|z||\big\}.
\end{eqnarray*}
We have
\begin{eqnarray}\nonumber
&&\frac{1}{2\pi}\int_{0}^{2\pi}\frac{1-|z|^2}{|e^{i\theta}-z|^2}\log\big(|f(e^{i\theta})-f(z/|z|)|+\delta\big)d\theta\\\nonumber&\leq&
\frac{1}{2\pi}\int_{\Gamma_1}
\frac{1-|z|^2}{|e^{i\theta}-z|^2}\log\big(\varepsilon\ \omega(|e^{i\theta}-z/|z||)+\delta\big)d\theta\\\nonumber
&&+\frac{1}{2\pi} \int_{\Gamma_2}
\frac{1-|z|^2}{|e^{i\theta}-z|^2}\log\big(\varepsilon\ \omega(|e^{i\theta}-z/|z||)|+\delta\big)d\theta\\\nonumber
&&+\frac{1}{2\pi} \int_{\Gamma_3}
\frac{1-|z|^2}{|e^{i\theta}-z|^2}\log\big(\omega(|e^{i\theta}-z/|z||)+\delta\big)d\theta\\\label{bo2}&:=& I_1+I_2+I_3.
\end{eqnarray}
It is clear that
\begin{equation}\label{bo3}
I_1\leq\big(\frac{1}{2\pi}\int_{\Gamma_1}\frac{1-|z|^2}{|e^{i\theta}-z|^2}d\theta\big) \log \big(\varepsilon\ \omega(1-|z|)+\delta\big).
\end{equation}
Next we have
\begin{eqnarray}\nonumber
    I_2&=&
\frac{1}{2\pi}\int_{\Gamma_2}\frac{1-|z|^2}{|e^{i\theta}-z|^2}\log\big(\varepsilon\ \omega(|e^{i\theta}-z/|z||)+\delta\big)d\theta\\\nonumber&\leq&
    \frac{1}{2\pi}\int_{\Gamma_2}\frac{1-|z|^2}{|e^{i\theta}-z|^2}\log\big(\varepsilon|e^{i\theta}-z/|z||\frac{\omega(1-|z|)}{1-|z|}+\delta\big)d\theta\\\nonumber&\leq&
\frac{1}{2\pi}\int_{\Gamma_2}\frac{1-|z|^2}{|e^{i\theta}-z|^2}\log\Big(\frac{|e^{i\theta}-z/|z||}{1-|z|}\big(\varepsilon\ \omega(1-|z|)+\delta\big)\Big)d\theta\\\nonumber&\leq&
\big(\frac{1}{2\pi}\int_{\Gamma_2}\frac{1-|z|^2}{|e^{i\theta}-z|^2}d\theta\big)\log\big(\varepsilon\ \omega(1-|z|)+\delta\big)\\\label{bo4}&&
+c \int_{t\geq 1}\frac{\log(t)}{t^2}dt,
\end{eqnarray}
where $c$ is an absolute constant.
Let $c^{'}_\varepsilon$ be a positif number such that $c^{'}_\varepsilon\leq c_\varepsilon$ and  for every $z\in\D$ satisfying $1-|z|\leq c^{'}_\varepsilon,$ we have
$\displaystyle \frac{1}{2\pi}\int_{\Gamma_3}\frac{1-|z|^2}{|e^{i\theta}-z|^2}d\theta\leq\varepsilon.$
Hence, as in \eqref{bo4}, we obtain
\begin{eqnarray}\nonumber
    I_3&=&
    \frac{1}{2\pi} \int_{\Gamma_3}
\frac{1-|z|^2}{|e^{i\theta}-z|^2}\log\big(\omega(|e^{i\theta}-z/|z||)+\delta\big)d\theta\\\nonumber&\leq&
\big(\frac{1}{2\pi}\int_{\Gamma_3}\frac{1-|z|^2}{|e^{i\theta}-z|^2}d\theta\big)\log\big(\omega(1-|z|)+\delta\big)+c \int_{t\geq 1}\frac{\log(t)}{t^2}dt\\\nonumber&=&
\big(\frac{1}{2\pi}\int_{\Gamma_3}\frac{1-|z|^2}{|e^{i\theta}-z|^2}d\theta\big)\log\big(\varepsilon\ \omega(1-|z|)+\delta\big)+ c \int_{t\geq 1}\frac{\log(t)}{t^2}dt \\\nonumber&&
+\big(\frac{1}{2\pi}\int_{\Gamma_3}\frac{1-|z|^2}{|e^{i\theta}-z|^2}d\theta\big)\log\Big(\frac{\omega(1-|z|)+\delta}{\varepsilon\ \omega(1-|z|)+\delta}\Big)
\\\nonumber&\leq&
\big(\frac{1}{2\pi}\int_{\Gamma_3}\frac{1-|z|^2}{|e^{i\theta}-z|^2}d\theta\big)\log\big(\varepsilon\ \omega(1-|z|)+\delta\big)+c \int_{t\geq 1}\frac{\log(t)}{t^2}dt
\\\label{bo5}&&
-\varepsilon\log(\varepsilon),
\end{eqnarray}
for every $z\in\D$ satisfying $1-|z|\leq c^{'}_\varepsilon.$
From  \eqref{bo3}, \eqref{bo4} and \eqref{bo5} we obtain
\begin{eqnarray}\nonumber
&&\exp\Big\{\frac{1}{2\pi}\int_{0}^{2\pi}\frac{1-|z|^2}{|e^{i\theta}-z|^2}\log\big(|f(e^{i\theta})-f(z/|z|)|+\delta\big)d\theta\Big\}
\\\label{bo6}&\leq&
A\big(\varepsilon\ \omega(|1-|z|)+\delta\big)\qquad(z\in\D \mbox{ and }1-|z|\leq c^{'}_\varepsilon).
\end{eqnarray}
where $A$ is an absolute constant. This completes the proof of the lemma.

\end{proof}

We use lemma \ref{key0} to prove the following one

\begin{lem}\label{Key}
Let $\omega$ be any arbitrary modulus of continuity and let $f\in\cL(\T).$   Then
\begin{equation}|f(z)-f(\zeta)|= o\big(\omega(|z-\zeta|)\big),\qquad(\mbox{as }|z-\zeta|\longrightarrow 0,\ z\in\D\text{ and }\zeta\in\T).
\end{equation}

\end{lem}

\begin{proof} We can suppose that $\|f\|_{\cL(\T)}\leq 1.$ Let $0<\varepsilon<1$ and $c_\varepsilon>0$ be a number such that for any $\xi,\zeta\in\T$ satisfying $|\xi-\zeta|\leq c_\varepsilon,$ we have
$|f(\xi)-f(\zeta)|\leq \varepsilon\ \omega(|\xi-\zeta|) .$ Fix $\zeta\in\T$ and fix $z\in\D$ such that $|z-\zeta|\leq c_\varepsilon/2$ and
$|z|\geq 1/4.$ We have $|z-\zeta|^2=(1-|z|)^2+|z||z/|z|-\zeta|^2\geq \frac{1}{4}|z/|z|-\zeta|^2.$
Hence $|z/|z|-\zeta|\leq 2|z-\zeta|\leq c_{\varepsilon}.$ We obtain
\begin{eqnarray*}\nonumber
&&|f(z)-f(\zeta)|\\\nonumber&\leq&|f(z)-f(z/|z|)|+|f(z/|z|)-f(\zeta)|\\\nonumber&\leq&
|f(z)-f(z/|z|)|+ \varepsilon\ \omega(|z/|z|-\zeta|)\\\label{bo1}&\leq&
\exp\Big\{\frac{1}{2\pi}\int_{0}^{2\pi}\frac{1-|z|^2}{|e^{i\theta}-z|^2}\log\big(|f(e^{i\theta})-f(z/|z|)|\big)d\theta\Big\}+ 2\varepsilon\ \omega(|z-\zeta|).
\end{eqnarray*}
Now, we use Lemma \ref{key0} to deduce the result of the lemma.
\end{proof}

\textit{\bf Proof of Theorem \ref{Tam}: }
Let $0<\varepsilon<1.$  From Lemma \ref{Key}
there is $0<c_\varepsilon<1/2$ such that for every $z\in\overline{\D}$ and for every $\zeta\in\T$
satisfying $|z-\zeta|\leq c_\varepsilon$ we have $|f(z)-f(\zeta)|\leq \varepsilon\ \omega(|z-\zeta|).$
Let $z,w\in\D$ satisfying  $|z-w|\leq c_{\varepsilon}/2$ and $\inf\{|z|,|w|\}=|w|\geq 1-c_\varepsilon.$ It follows that $|z/|z|-w/|w||\leq |zw|^{-1/2}|z-w|\leq2|z-w|\leq c_{\varepsilon}.$
\begin{itemize}
  \item [$1.$] First we assume that $|z-w|\geq 1-|w|.$
We obtain
\begin{eqnarray}\nonumber
&&|f(z)-f(w)|\\\nonumber&\leq& |f(z)-f(z/|z|)|+|f(z/|z|)-f(w/|w|)|+|f(w)-f(w/|w|)|\\\label{T1} &\leq&
4\varepsilon\ \omega(|z-w|).
\end{eqnarray}
  \item [$2.$] Now we suppose that $|z-w|\leq 1-|w|.$
We apply the maximum principle theorem in $\D$ to the analytic function
$\displaystyle z\mapsto \frac{f(z)-f(w)}{z-w},$  we get
\begin{eqnarray*}
\Big|\frac{f(z)-f(w)}{z-w}\Big|\leq\sup_{\xi\in\T}\Big|\frac{f(\xi)-f(w)}{\xi-w}\Big|=
\frac{|f(\xi_w)-f(w)|}{|\xi_w-w|} \qquad(\xi_w\in\T).
\end{eqnarray*}
\begin{itemize}
  \item [$2.1.$] If $|\xi_w-w|\geq c_\varepsilon,$ then
\begin{eqnarray}\nonumber
\frac{|f(z)-f(w)|}{\omega(|z-w|)}&=&\nonumber
\Big|\frac{f(z)-f(w)}{z-w}\Big|\frac{|z-w|}{\omega(|z-w|)}\\\nonumber&\leq&
2\frac{\|f\|_{\infty}}{c_\varepsilon} \frac{|z-w|}{\omega(|z-w|)}\\\label{T2}&=&
o(1)\qquad(\mbox{as }|z-w|\longrightarrow0).
\end{eqnarray}
  \item [$2.2.$]If $|\xi_w-w|\leq c_\varepsilon,$ then $|f(\xi_w)-f(w)|\leq \varepsilon\ \omega(|\xi_w-w|).$ It follows that
\begin{eqnarray}\nonumber
\frac{|f(z)-f(w)|}{\omega(|z-w|)}&=&
\Big|\frac{f(z)-f(w)}{z-w}\Big|\frac{|z-w|}{\omega(|z-w|)}\\\nonumber&\leq&
\Big|\frac{f(\xi_w)-f(w)}{\xi_w-w}\Big|\frac{1-|w|}{\omega(1-|w|)}\qquad(\xi_w\in\T)\\\nonumber&\leq&
\varepsilon\frac{\omega(|\xi_w-w|)}{|\xi_w-w|}\frac{1-|w|}{\omega(1-|w|)}\\\label{T3}&\leq&
\varepsilon.
\end{eqnarray}
\end{itemize}
From inequalities \eqref{T1}, \eqref{T2} and \eqref{T3}, there is $c^{'}_\varepsilon>0$ such that if $|z-w|\leq c^{'}_\varepsilon$ and  $\inf\{|z|,|w|\}\geq 1-c^{'}_\varepsilon,$ then
\begin{eqnarray}\label{1step}
\frac{|f(z)-f(w)|}{\omega(|z-w|)}&\leq& \varepsilon.
\end{eqnarray}
\item [$3.$] If $z,w\in\D$ are such that $\sup\{|z|,|w|\}\leq 1-c^{'}_\varepsilon$ we have
\begin{eqnarray}\nonumber
\frac{|f(z)-f(w)|}{\omega(|z-w|)}&=&\Big|\frac{f(z)-f(w)}{z-w}\Big|\frac{|z-w|}{\omega(|z-w|)}\\\nonumber
&\leq& \sup_{|\zeta|\leq 1-c^{'}_\varepsilon}|f'(\zeta)|\frac{|z-w|}{\omega(|z-w|)}\\\label{2step}
&=& o(1)\qquad(\mbox{as }|z-w|\longrightarrow0).
\end{eqnarray}
\end{itemize}
From \eqref{1step} and \eqref{2step} we deduce the result. The proof of the theorem is completed.

\section{\bf Appendix B. Factorization property in $\cL.$}

The F-property of $\cL,$ for any arbitrary modulus of continuity $\omega,$  is given by Shirokov \cite{Shi2}. For completeness we give here the proof

\begin{thm}\cite{Shi2}\label{fpr}
Let $\omega$ be any arbitrary modulus of continuity. Let $f$ be a function in $\cL$ and let $U$ be inner function such that $f/U\in\cH^{\infty}(\D).$ Then $f/U\in\cL$ and
$$\frac{|f/U(z)-f/U(w)|}{\omega(|z-w|)}=o(1)\qquad(as\ |z-w|\longrightarrow 0),$$
uniformly with respect to $U.$
Also $\|f/U\|_{\omega}\leq c\|f\|_{\omega},$ where $c$ is an absolute constant.
\end{thm}

\begin{coro}\label{fprc}
Let $\omega$ be any arbitrary modulus of continuity. Let $f$ and $g$ be functions in $\cL$ and let $\{U_n\}_{n\in\N}$ be sequence of inner functions such that $f/U_n\in\cH^{\infty}(\D)$ for every $n\in\N.$ If $\lim\limits_{n\rightarrow +\infty}\|f/U_n-g\|_{\infty}=0,$ then
$\lim\limits_{n\rightarrow +\infty}\|f/U_n-g\|_{\omega}=0.$
\end{coro}

\begin{proof}
The proof immediately follows from Theorem \ref{fpr} and Lemma
\ref{shilbert}.
\end{proof}

We begin the proof of Theorem \ref{fpr} by establishing several lemmas.

\begin{lem}\label{Fpr1}
Let $\omega$ be any arbitrary modulus of continuity. Let $f$ be a function in $\cL.$  Then
$$|O_{_f}(z)|\leq o(\omega(1-|z|))+A|f(z/|z|)|\qquad(\mbox{as }|z|\longrightarrow 1),$$
where $A>0$ is an absolute constant.
\end{lem}

\begin{proof}
For $z\in\D,$ we have
\begin{eqnarray*}
&&\log|O_{_f}(z)|
\\&=&\frac{1}{2\pi}\int_{0}^{2\pi}\frac{1-|z|^2}{|e^{i\theta}-z|^2}\log|f(e^{i\theta})|d\theta\\&\leq&
\frac{1}{2\pi}\int_{0}^{2\pi}\frac{1-|z|^2}{|e^{i\theta}-z|^2}\log\big(|f(e^{i\theta})-f(z/|z|)|+|f(z/|z|)|\big)d\theta.
\end{eqnarray*}
Now, we use Lemma \ref{key0} to complete the proof of the lemma.
\end{proof}

For a function $f\in\cH^{\infty}(\D)$ we set
$$ a_{_f}(\xi):=\sum\limits_{n\geq0} \frac{1-|a_n|^2}{|\xi- a_n|^2}+
\frac{1}{\pi}\int_{\T}\frac{1}{|e^{i\theta}-\xi|^2}d\mu_{_f}(\theta),$$
where $\{a_n\ :\ n\in\N\}=Z_{_f}\cap\D$ (for all $n,$ $a_n$ is repeated according to it's multiplicity) and $\mu_{_f}$ is the positive singular measure associated to the singular factor $S_{_f}$ of $f.$

\begin{lem}\label{Fpr2}
Let $f$ be a function in the disk algebra with inner factor $U_{_f}\not\equiv1.$   Let $\xi\in\T\setminus E_{_f}$ and let $0<\rho<1$ be such that  $1-\rho\leq d(\xi, Z_{_f}).$ Then
$$|U_{_f}(\rho \xi)|\leq \exp\big\{-\frac{1-\rho}{8}a_{_f}(\xi)\big\}.$$
\end{lem}

\begin{proof}
We have
\begin{eqnarray}\nonumber
\log|S_{_f}(\rho \xi )|&=&-\frac{1}{2\pi}\int_{\T}\frac{1-|\rho \xi|^2}{|e^{i\theta}-\rho \xi|^2}d\mu_{_f}(\theta)\\\nonumber&\leq&
-\frac{1}{2\pi}\int_{\T}\frac{1}{4}\frac{1-\rho}{|e^{i\theta}-\xi|^2}d\mu_{_f}(\theta)\\\label{esS}&=&
-\frac{1-\rho}{8}\frac{1}{\pi}\int_{\T}\frac{1}{|e^{i\theta}-\xi|^2}d\mu_{_f}(\theta).
\end{eqnarray}
It is clear that $$|z-w|^2=||z|-|w||^2+|z||w||z/|z|-w/|w||^2 \qquad(z,w\in\D).$$
Now we estimate $|B_{_f}(\rho \xi)|.$ For all $n\in\N,$ we have
\begin{eqnarray*}
    \Big|\frac{\rho \xi- a_n}{\xi-\rho a_n}\Big|^2&=&\frac{(\rho-|a_n|)^2+\rho|a_n||\xi-a_n/|a_n||^2}{(1-\rho|a_n|)^2+\rho|a_n||\xi-a_n/|a_n||^2}\\&=&
1-(1-\rho^2)\frac{1-|a_n|^2}{|\xi-\rho a_n|^2}\\&\leq&
1-\frac{1-\rho}{4}\frac{1-|a_n|^2}{|\xi- a_n|^2}.
\end{eqnarray*}
Hence $$\log \Big|\frac{\rho \xi- a_n}{\xi-\rho a_n}\Big|\leq -\frac{1-\rho}{8}\frac{1-|a_n|^2}{|\xi- a_n|^2}.
$$
Therefore
\begin{equation}\label{esB}|B_{_f}(\rho \xi)|\leq \exp\big\{-\frac{1-\rho}{8}\sum\limits_{n\geq0} \frac{1-|a_n|^2}{|\xi- a_n|^2}\big\}.
\end{equation}
From \eqref{esS} and \eqref{esB} we obtain
\begin{equation*}
|U_{_f}(\rho \xi)|\leq \exp\big\{-\frac{1-\rho}{8}a_{_f}(\xi)\big\}.
\end{equation*}
This proves the lemma.
\end{proof}

\begin{lem}\label{derni�re}
Let $f$ be a function in $\cL$ and let $U$ be an inner function such that $f/U\in\cH^{\infty}(\D).$ Then
$$|f(\zeta)|=o\big(\omega(\frac{1}{a_{_U}(\zeta)})\big)\qquad(\mbox{as }d(\zeta,Z_{_f})\longrightarrow0, \ \zeta\in\T\setminus E_{_f}).$$
\end{lem}

\begin{proof}
Let $\varepsilon>0.$ There is $c_\varepsilon>0$ such that if $|z-w|\leq c_\varepsilon,$ $z,w\in\T,$ then we have $|f(z)-f(w)|\leq \varepsilon\omega(|z-w|).$ From Lemma \ref{Fpr1}, there is  $0<c^{'}_\varepsilon<c_\varepsilon$ such that if $1-\rho\leq c^{'}_\varepsilon,$ then $|O_f(\rho \zeta)|\leq A\big(\varepsilon\ \omega(1-\rho)+|f(\zeta)|\big),$
where $A\geq1$ is an absolute constant. Let $\zeta\in\T\setminus E_{_f}$ such that $d(\zeta, Z_{_f})\leq c^{'}_\varepsilon.$
\begin{enumerate}
  \item [1.] We assume that $\displaystyle a_{_U}(\zeta)\leq \frac{8A}{d(\zeta, Z_{_f})}.$ Then we obtain
 $\displaystyle|f(\zeta)|\leq \varepsilon\omega(d(\zeta, Z_{_f}))\leq\varepsilon \omega(\frac{8A}{a_{_U}(\zeta)})\leq 8A\varepsilon\omega(\frac{1}{a_{_U}(\zeta)}).$
  \item [2.] Now assume that $\displaystyle a_{_U}(\zeta)\geq \frac{8A}{d(\zeta, Z_{_f})}.$
 Then $1-\rho_{_\zeta}\leq d(\zeta, Z_{_f}),$ where $\displaystyle \rho_{_\zeta}:=1-\frac{8A}{a_{_U}(\zeta)}.$
 We have $U_{_f}/U\in\cH^{\infty}(\D),$ then $\displaystyle a_{_U}(\zeta)\leq a_{_f}(\zeta)$ and by using Lemma \ref{Fpr2} we obtain
$$|U_{_f}(\rho_{_\zeta} \zeta)|\leq\exp\big\{-\frac{1-\rho_{_\zeta}}{8}a_{_f}(\zeta)\big\}\leq\exp\big\{-\frac{1-\rho_{_\zeta}}{8}a_{_U}(\zeta)\big\} =\exp\{-A\}.$$ Since $d(\zeta, Z_{_f})\leq c^{'}_\varepsilon,$ then  $1-\rho_{_\zeta}\leq c^{'}_\varepsilon$ and we have
$$|f(\rho_{_\zeta} \zeta)|=|U_f(\rho_{_\zeta} \zeta)||O_f(\rho_{_\zeta} \zeta)|\leq A\exp\{-A\}\big(\varepsilon\ \omega(1-\rho_{_\zeta})+|f(\zeta)|\big).$$
Hence
\begin{eqnarray*} |f(\zeta)|&\leq& |f(\zeta)-f(\rho_{_\zeta} \zeta)|+|f(\rho_{_\zeta} \zeta)|\\&\leq&
\varepsilon\ \omega(1-\rho_{_\zeta})+ A\exp\{-A\}\big(\varepsilon\ \omega(1-\rho_{_\zeta})+|f(\zeta)|\big).
\end{eqnarray*}
It follows that
$\displaystyle|f(\zeta)|\leq 3\varepsilon\ \omega(1-\rho_{_\zeta})\leq 24 A\varepsilon\ \omega(\frac{1}{a_{_U}(\zeta)}).$ This completes the proof of the lemma.
\end{enumerate}
\end{proof}

\textit{ \bf Proof of Theorem \ref{fpr} :}

Now, we
can deduce the proof of Theorem \ref{fpr} by using
Lemma \ref{derni�re}. Indeed, from Theorem \ref{Tam}
it is sufficient to prove that $f/U\in\cL(\T),$ that is
\begin{equation*}
\frac{|f(\xi)/U(\xi)-f(\zeta)/U(\zeta)|}{\omega(|\xi-\zeta|)}=o(1)\qquad(\text{as } |\xi-\zeta|\longrightarrow 0\mbox{ and }\xi,\zeta\in\T ).
\end{equation*}
Let $\xi,\zeta\in\T$ be two distinct points such that $d(\xi, Z_{_f})\geq d(\zeta, Z_{_f}).$ We have
\begin{eqnarray*}
\frac{|f(\xi)/U(\xi)-f(\zeta)/U(\zeta)|}{\omega(|\xi-\zeta|)} &\leq&
\frac{|f(\xi)-f(\zeta)|}{\omega(|\xi-\zeta|)}+ |f(\zeta)|\frac{|U(\xi)-U(\zeta)|}{\omega(|\xi-\zeta|)}.
\end{eqnarray*}
Then it suffices to prove
\begin{equation}\label{fprbp}
|f(\zeta)|\frac{|U(\xi)-U(\zeta)|}{\omega(|\xi-\zeta|)}=o(1)\qquad(\mbox{as }|\xi-\zeta|\longrightarrow0).
\end{equation}

\begin{enumerate}
\item [$1.$] First we suppose that $|\xi-\zeta|\geq \frac{1}{2}d(\zeta, Z_{_f}).$ Then
\begin{eqnarray}\nonumber
|f(\zeta)|\frac{|U(\xi)-U(\zeta)|}{\omega(|\xi-\zeta|)} &\leq&
2\frac{|f(\zeta)|}{\omega(\frac{1}{2}d(\zeta, Z_{_f}))} \\\nonumber&\leq&
4\frac{|f(\zeta)|}{\omega(d(\zeta, Z_{_f}))}\\\label{fprb1}
&=&o(1)\qquad (\mbox{as}\ |\xi-\zeta|\longrightarrow 0).
\end{eqnarray}
\item [$2.$] Next we suppose that $|\xi-\zeta|\leq \frac{1}{2}d(\zeta, Z_{_f}).$ Then $[\xi,\zeta]\subset
T\setminus E_{_f}.$ There is $z\in[\xi,\zeta]$ such that $\displaystyle \frac{|U(\xi)-U(\zeta)|}{|\xi-\zeta|}=|U'(z)|\leq\sum\limits_{n\geq0}\frac{1-|a_n|^2}{|z-a_n|^2}
+\frac{1}{\pi}\int_{\T}\frac{1}{|e^{i\theta}-z|^2}d\mu_{_f}(\theta):=a_{_U}(z)\leq 4a_{_U}(\zeta)
\leq \frac{c_{_f}}{d^{2}(\zeta, Z_{_f})}.$
\begin{itemize}
\item [$2.1.$] If  $\displaystyle a_{_U}(\zeta)\leq \frac{1}{|\xi-\zeta|}.$ Then, by using Lemma \ref{derni�re}, we deduce that
$$|f(\zeta)|a_{_U}(\zeta)\frac{|\xi-\zeta|}{\omega(|\xi-\zeta|)}=
o(1)\qquad(\mbox{as } d(\zeta, Z_{_f})\longrightarrow0).$$
Therefore
\begin{eqnarray}\nonumber
&&|f(\zeta)|\frac{|U(\xi)-U(\zeta)|}{\omega(|\xi-\zeta|)}\\\nonumber&\leq&
4|f(\zeta)|a_{_U}(\zeta)\frac{|\xi-\zeta|}{\omega(|\xi-\zeta|)}\\\nonumber&\leq&
\inf\Big\{4|f(\zeta)|a_{_U}(\zeta)\frac{|\xi-\zeta|}{\omega(|\xi-\zeta|)},\ \frac{c_{_f}\|f\|_{\infty}}{d^2(\zeta, Z_{_f})}\frac{|\xi-\zeta|}{\omega(|\xi-\zeta|)}\Big\}\\\label{fprb2}&=&
o(1)\qquad(\mbox{as } |\xi-\zeta|\longrightarrow0).
\end{eqnarray}
\item [$2.2.$] Now assume that $\displaystyle a_{_U}(\zeta)\geq \frac{1}{|\xi-\zeta|}.$ Then, by  Lemma \ref{derni�re}, we obtain
$|f(\zeta)|=o(\omega(|\xi-\zeta|)),$ as $d(\zeta, Z_{_f})\longrightarrow0.$ Therefore
\begin{eqnarray}\nonumber
|f(\zeta)|\frac{|U(\xi)-U(\zeta)|}{\omega(|\xi-\zeta|)}&\leq&
\inf\Big\{2\frac{|f(\zeta)|}{\omega(|\xi-\zeta|)},\ \frac{c_{_f}\|f\|_{\infty}}{d^2(\zeta, Z_{_f})}\frac{|\xi-\zeta|}{\omega(|\xi-\zeta|)}\Big\}\\\label{fprb4}&=&
o(1)\qquad(\mbox{as }|\xi-\zeta|\longrightarrow0).
\end{eqnarray}
\end{itemize}
\end{enumerate}
Consequently \eqref{fprbp} follows from inequalities \eqref{fprb1},
\eqref{fprb2} and  \eqref{fprb4}. The proof of the theorem is completed.

\vspace*{1cm}

\textsc{Acknowledgements.} I wish to thank Professors A. Borichev,
O. El-Fallah and K. Kellay   for the great interest which they carried to
this work.

\vspace*{1cm}

\end{document}